\documentclass[12pt]{article}
\usepackage{fullpage,amssymb,amsmath,amsthm}
\usepackage{ifpdf}
\ifpdf
 \usepackage[hyperindex,pagebackref]{hyperref}%
\else
 \expandafter\ifx\csname dvipdfm\endcsname\relax
 \usepackage[hypertex,hyperindex,pagebackref]{hyperref}
 \else
 \usepackage[dvipdfm,hyperindex,pagebackref]{hyperref}
 \fi
\fi
\allowdisplaybreaks[4]
\theoremstyle{plain}
\newtheorem{thm}{Theorem}
\theoremstyle{remark}

\DeclareMathOperator{\td}{d}
\newcommand{\bell}{\textup{B}}

\begin{document}

\title{An Explicit Formula for Bell Numbers in Terms of Stirling Numbers and Hypergeometric Functions}

\author{Feng Qi}

\date{Department of Mathematics, College of Science\\ Tianjin Polytechnic University\\ Tianjin City, 300387, \\ People's Republic of China\\
\href{mailto: F. Qi <qifeng618@gmail.com>}{qifeng618@gmail.com}\\ \href{mailto: F. Qi <qifeng618@hotmail.com>}{qifeng618@hotmail.com}\\ \href{mailto: F. Qi <qifeng618@qq.com>}{qifeng618@qq.com}\\
\url{http://qifeng618.wordpress.com}}

\maketitle

\begin{abstract}
In the paper, the author finds an explicit formula for computing Bell numbers in terms of Kummer confluent hypergeometric functions and Stirling numbers of the second kind.
\end{abstract}

\section{Introduction}

In combinatorics, Bell numbers, usually denoted by $\bell_n$ for $n\in\{0\}\cup\mathbb{N}$, count the number of ways a set with $n$ elements can be partitioned into disjoint and non-empty subsets. These numbers have been studied by mathematicians since the $19$th century, and their roots go back to medieval Japan, but they are named after Eric Temple Bell, who wrote about them in the $1930$s.
Every Bell number $\bell_n$ may be generated by
\begin{equation}\label{Bell-generate-function}
e^{e^x-1}=\sum_{k=0}^\infty \bell_k\frac{x^k}{k!}
\end{equation}
or, equivalently, by
\begin{equation}\label{Bell-generate-funct-2nd}
e^{e^{-x}-1}=\sum_{k=0}^\infty(-1)^k\bell_k\frac{x^k}{k!}.
\end{equation}
\par
In combinatorics, Stirling numbers arise in a variety of combinatorics problems. They are introduced in the eighteen century by James Stirling. There are two kinds of Stirling numbers: Stirling numbers of the first kind and Stirling numbers of the second kind.
Every Stirling number of the second kind, usually denoted by $S(n,k)$, is the number of ways of partitioning a set of $n$ elements into $k$ nonempty subsets, may be computed by
\begin{equation}\label{S(n-k)-explicit}
S(n,k)=\frac1{k!}\sum_{i=0}^k(-1)^i\binom{k}{i}(k-i)^n,
\end{equation}
and may be generated by
\begin{equation}\label{2stirling-gen-funct-exp}
\frac{(e^x-1)^k}{k!}=\sum_{n=k}^\infty S(n,k)\frac{x^n}{n!}, \quad k\in\{0\}\cup\mathbb{N}.
\end{equation}
\par
In the theory of special functions, the hypergeometric functions are denoted and defined by
\begin{equation}\label{hypergeom-f}
{}_pF_q(a_1,\dotsc,a_p;b_1,\dotsc,b_q;x)=\sum_{n=0}^\infty\frac{(a_1)_n\dotsm(a_p)_n} {(b_1)_n\dotsm(b_q)_n}\frac{x^n}{n!}
\end{equation}
for $b_i\notin\{0,-1,-2,\dotsc\}$ and $p,q\in\mathbb{N}$, where $(a)_0=1$ and $(a)_n=a(a+1)\dotsm(a+n-1)$ for $n\in\mathbb{N}$ and any complex number $a$ is called the rising factorial. Specially, the series
\begin{equation}
{}_1F_1(a;b;z)=\sum_{k=0}^\infty\frac{(a)_k}{(b)_k}\frac{z^k}{k!}
\end{equation}
is called Kummer confluent hypergeometric function.
\par
In combinatorics or number theory, it is common knowledge that Bell numbers $\bell_n$ may be computed in terms of Stirling numbers of the second kind $S(n,k)$ by
\begin{equation}\label{Bell-Stirling-eq}
\bell_n=\sum_{k=1}^nS(n,k).
\end{equation}
\par
In this paper, we will find a new explicit formula for computing Bell numbers $\bell_n$ in terms of Kummer confluent hypergeometric functions ${}_1F_1(k+1;2;1)$ and Stirling numbers of the second kind $S(n,k)$ as follows.

\begin{thm}\label{Bell-Stirling-HyperG-thm}
For $n\in\mathbb{N}$, Bell numbers $\bell_n$ may be expressed as
\begin{equation}
\bell_n =\frac1e\sum_{k=1}^n(-1)^{n-k}S(n,k)k!{}_1F_1(k+1;2;1).
\end{equation}
\end{thm}

\section{Proof of Theorem~\ref{Bell-Stirling-HyperG-thm}}

We now start out to verify Theorem~\ref{Bell-Stirling-HyperG-thm} as follows.
\par
Among other things, Qi and Wang~\cite[Theorem~1.2]{simp-exp-degree-revised.tex} obtained that the function
\begin{equation}\label{exp=k=sum-eq-degree=k+1}
H_k(z)=e^{1/z}-\sum_{m=0}^k\frac{1}{m!}\frac1{z^m}
\end{equation}
for $k\in\{0\}\cup\mathbb{N}$ and $z\ne0$ has the integral representation
\begin{equation}\label{exp=k=degree=k+1-int}
H_k(z)=\frac1{k!(k+1)!}\int_0^\infty {}_1F_2(1;k+1,k+2;t)t^k e^{-zt}\td t, \quad \Re(z)>0.
\end{equation}
See also~\cite[Section~1.2]{Bessel-ineq-Dgree-CM.tex} and~\cite[Lemma~2.1]{QiBerg.tex}. When $k=0$, the integral representation~\eqref{exp=k=degree=k+1-int} becomes
\begin{equation}\label{open-answer-1}
e^{1/z}=1+\int_0^\infty \frac{I_1\bigl(2\sqrt{t}\,\bigr)}{\sqrt{t}\,} e^{-zt}\td t, \quad \Re(z)>0,
\end{equation}
where $I_\nu(z)$ stands for the modified Bessel function of the first kind
\begin{equation}\label{I=nu(z)-eq}
I_\nu(z)= \sum_{k=0}^\infty\frac1{k!\Gamma(\nu+k+1)}\biggl(\frac{z}2\biggr)^{2k+\nu}
\end{equation}
for $\nu\in\mathbb{R}$ and $z\in\mathbb{C}$, see~\cite[p.~375, 9.6.10]{abram}, and $\Gamma$ represents the classical Euler gamma function which may be defined by
\begin{equation}\label{gamma-dfn}
\Gamma(z)=\int^\infty_0t^{z-1} e^{-t}\td t, \quad \Re z>0,
\end{equation}
see~\cite[p.~255]{abram}. Replacing $z$ by $e^x$ in~\eqref{open-answer-1} gives
\begin{equation}\label{open-answer-exp}
e^{1/e^{x}}=e^{e^{-x}}=1+\int_0^\infty \frac{I_1\bigl(2\sqrt{t}\,\bigr)}{\sqrt{t}\,} e^{-e^{x}t}\td t.
\end{equation}
Differentiating $n\ge1$ times with respect to $x$ on both sides of~\eqref{open-answer-exp} and~\eqref{Bell-generate-funct-2nd} gives
\begin{equation}\label{open-answer-exp-diff}
\frac{\td^ne^{e^{-x}}}{\td x^n}=\int_0^\infty \frac{I_1\bigl(2\sqrt{t}\,\bigr)}{\sqrt{t}\,} \frac{\td^ne^{-e^{x}t}}{\td x^n}\td t
\end{equation}
and
\begin{equation}\label{Bell-generate-funct-diff}
\frac{\td^ne^{e^{-x}}}{\td x^n}=e\sum_{k=n}^\infty(-1)^k\bell_k\frac{x^{k-n}}{(k-n)!}.
\end{equation}
 From~\eqref{open-answer-exp-diff} and~\eqref{Bell-generate-funct-diff}, it follows that
\begin{equation*}
e\sum_{k=n}^\infty(-1)^k\bell_k\frac{x^{k-n}}{(k-n)!}
=\int_0^\infty \frac{I_1\bigl(2\sqrt{t}\,\bigr)}{\sqrt{t}\,} \frac{\td^ne^{-e^{x}t}}{\td x^n}\td t.
\end{equation*}
Taking $x\to0$ in the above equation yields
\begin{equation}\label{bell-int-rep-eq}
(-1)^ne\bell_n =\int_0^\infty \frac{I_1\bigl(2\sqrt{t}\,\bigr)} {\sqrt{t}\,} \lim_{x\to0}\frac{\td^ne^{-e^{x}t}}{\td x^n}\td t.
\end{equation}
\par
In combinatorics, Bell polynomials of the second kind, or say, the partial Bell polynomials, $\bell_{n,k}(x_1,x_2,\dotsc,x_{n-k+1})$ are defined by
\begin{equation}
\bell_{n,k}(x_1,x_2,\dotsc,x_{n-k+1})=\sum_{\substack{1\le i\le n,\ell_i\in\mathbb{N}\\ \sum_{i=1}^ni\ell_i=n\\ \sum_{i=1}^n\ell_i=k}}\frac{n!}{\prod_{i=1}^{n-k+1}\ell_i!} \prod_{i=1}^{n-k+1}\Bigl(\frac{x_i}{i!}\Bigr)^{\ell_i}
\end{equation}
for $n\ge k\ge1$, see~\cite[p.~134, Theorem~A]{Comtet-Combinatorics-74}, and satisfy
\begin{equation}\label{Bell(n-k)}
\bell_{n,k}\bigl(abx_1,ab^2x_2,\dotsc,ab^{n-k+1}x_{n-k+1}\bigr) =a^kb^n\bell_{n,k}(x_1,x_n,\dotsc,x_{n-k+1})
\end{equation}
and
\begin{equation}\label{Bell-stirling}
\bell_{n,k}(\overbrace{1,1,\dotsc,1}^{n-k+1})=S(n,k),
\end{equation}
see~\cite[p.~135]{Comtet-Combinatorics-74}, where $a$ and $b$ are any complex numbers. The well-known Fa\`a di Bruno formula may be described in terms of Bell polynomials of the second kind $\bell_{n,k}(x_1,x_2,\dotsc,x_{n-k+1})$ by
\begin{equation}\label{Bruno-Bell-Polynomial}
\frac{\td^n}{\td x^n}f\circ g(x)=\sum_{k=1}^nf^{(k)}(g(x)) \bell_{n,k}\bigl(g'(x),g''(x),\dotsc,g^{(n-k+1)}(x)\bigr),
\end{equation}
see~\cite[p.~139, Theorem~C]{Comtet-Combinatorics-74}.
By Fa\`a di Bruno formula~\eqref{Bruno-Bell-Polynomial} and the identities~\eqref{Bell(n-k)} and~\eqref{Bell-stirling}, we have
\begin{equation}\label{der-lin-fin-eq}
\begin{aligned}
\frac{\td^ne^{-e^{x}t}}{\td x^n}
&=\sum_{k=1}^n e^{-e^{x}t} \bell_{n,k}(\overbrace{-e^{x}t,-e^{x}t,\dotsc,-e^{x}t}^{n-k+1})\\
&= e^{-e^{x}t}\sum_{k=1}^n(-e^{x}t)^k \bell_{n,k}(\overbrace{1,1,\dotsc,1}^{n-k+1})\\
&= e^{-e^{x}t}\sum_{k=1}^n(-e^{x}t)^kS(n,k)\\
&\to e^{-t}\sum_{k=1}^n(-t)^kS(n,k)
\end{aligned}
\end{equation}
as $x\to0$. Substituting~\eqref{der-lin-fin-eq} into~\eqref{bell-int-rep-eq} leads to
\begin{align*}
\bell_n &=\frac1e\sum_{k=1}^n(-1)^{n-k}S(n,k)\int_0^\infty I_1\bigl(2\sqrt{t}\,\bigr)t^{k-1/2}e^{-t}\td t\\
&=\frac1e\sum_{k=1}^n(-1)^{n-k}S(n,k)k!{}_1F_1(k+1;2;1).
\end{align*}
The required proof is complete.

\bigskip
\hrule height1pt width\textwidth\bigskip\noindent
2010 \emph{Mathematics Subject Classification}: Primary 11B73; Secondary 33B10, 33C15.
\par\noindent
\emph{Keywords}: explicit formula, Bell number, confluent hypergeometric function of the first kind, Stirling number of the second kind.\bigskip
\hrule height1pt width\textwidth\bigskip\noindent
(Concerned with sequences \href{https://oeis.org/A000110}{A000110} and  \href{https://oeis.org/A008277}{A008277}.)\bigskip\noindent
\hrule height1pt width\textwidth\bigskip\noindent

\end{document}